\label{key}
\documentclass[english, 11pt]{article}

\usepackage[T1]{fontenc}
\usepackage[latin9]{inputenc}
\usepackage{amsthm}
\usepackage{amsmath}
\usepackage{amssymb}
\usepackage{color}
\usepackage{esint}
\usepackage{verbatim}

\makeatletter
\numberwithin{equation}{section}
\numberwithin{figure}{section}
\theoremstyle{plain}
\newtheorem{thm}{\protect\theoremname}
\theoremstyle{definition}
\newtheorem{defn}[thm]{\protect\definitionname}
\theoremstyle{remark}
\newtheorem{rem}[thm]{\protect\remarkname}
\theoremstyle{plain}

\@ifundefined{date}{}{\date{}}
\makeatother

\usepackage{babel}
  \providecommand{\corollaryname}{Corollary}
  \providecommand{\definitionname}{Definition}
\providecommand{\remarkname}{Remark}
\providecommand{\theoremname}{Theorem}

\def\pat{\partial_t}

\def\pay{\partial_y}
\def\vw{{\bf w}}
\def\vv{{\bf v}}
\def\vu{{\bf u}}
\def\vn{{\bf n}}

\def\vw{{\bf w}}

\def\vV{{\bf V}}
\def\vU{{\bf U}}
\def\vW{{\bf W}}

\def\vtau{\boldsymbol{\tau}}

\def\diver{{\operatorname{div}}}

\begin{document}

\title{Dimension reduction for the compressible Navier-Stokes system with density dependent viscosity}

\author{Matteo Caggio${^1}$, V\' aclav M\' acha${^2}$}
\date{${^1}$Department of Mathematics, Faculty of Science, University of Zagreb, Croatia \\ ${^2}$Institute of Mathematics, AV\v CR, Prague, Czech Republic \\ 
\smallskip \small matteo.caggio@math.hr, macha@math.cas.cz}

\maketitle
\begin{abstract}
	We consider a compressible Navier-Stokes system for a barotropic fluid with density dependent viscosity in a three-dimensional time-space domain $(0,T)\times \Omega_\varepsilon$ where $\Omega_\varepsilon = (0,\varepsilon)^2\times (0,1)$. We show that the weak solutions of the 3D system converges to the strong solution of the respective 1D system as $\varepsilon\rightarrow0.$ 
\end{abstract}
\smallskip


\tableofcontents{}

\newpage{}

\section{Introduction}

The present paper is devoted to the problem of the limit passage from three-dimensional to one-dimensional geometry. In a three-dimensional time-space domain $(0,T)\times \Omega_\varepsilon$, where 
\begin{equation} \label{space-dom}
\Omega_\varepsilon = (0,\varepsilon)^2\times (0,1), \ \ \varepsilon > 0,
\end{equation}
we consider the following compressible Navier-Stokes system with density dependent viscosity coefficients
\begin{equation} \label{cont}
\pat \varrho + \diver(\varrho \vu)  = 0,
\end{equation}
\begin{equation} \label{mom}
\pat(\varrho \vu) + \diver (\varrho\vu \otimes\vu) + \nabla p(\varrho) - \diver (2\mu(\varrho) D(\vu)) - \nabla (\lambda(\varrho) \diver \vu) = 0,
\end{equation}
supplemented with the initial conditions
\begin{equation} \label{ic}
    \varrho = \varrho (0,\cdot), \ \ \varrho\vu = \varrho\vu (0, \cdot)
\end{equation}
and the mixed boundary conditions 
\begin{equation}
\begin{split}
&\vu\restriction_{\partial_1\Omega}\cdot \vn = 0, \ \ (D(\vu)\restriction_{\partial_1\Omega}\vn)\cdot \vtau = 0,\\
&\vu(t,x_1,x_2,0) = \vu(t,x_1,x_2,1), \ \  \varrho(t,x_1,x_2,0) = \varrho(t,x_1,x_2,1),
\end{split}
\end{equation}
where $\vn$ is outer unit normal vector, $\vtau$ denotes an arbitrary tangential vector and
\begin{equation}
\partial_1\Omega = ((0,\varepsilon)\times \{0,\varepsilon\}\times (0,1)) \cup (\{0,\varepsilon\}\times(0,\varepsilon)\times(0,1)).
\end{equation}
and the second boundary condition holds for all $t\in(0,T)$ and all $(x_1,x_2)\in (0,\varepsilon)^2$.

The above boundary conditions are a mix of the full-slip boundary conditions on the "long walls" and the periodic boundary conditions on the "shrinking walls". The full-slip boundary conditions are the only relevant for the dimension reduction as the use of no-slip boundary conditions gives only trivial result. On the other hand, the periodic boundary condition helps to handle certain boundary terms involving pressure (see Section 3, relation (3.20)).

Here, $\varrho=\varrho(t,x)$ and $\vu=\vu(t,x)$ stand for the density and the velocity field of the fluid. The quantity $D(\vu)=(\nabla\vu+\nabla^t \vu)/2$ represents the strain tensor. We denote by $\mu(\varrho)$ and $\lambda(\varrho)$ the viscosity coefficients functions of the fluid density. We consider a pressure $p$ of the type $p(\varrho)=a\varrho^\gamma$ where $a>0$ and $\gamma>1$.

In the context of thin-domains limit analysis for compressible fluids with density dependent viscosity coefficients, a recent paper of Zhang \cite{Z} provides a reduction from two to one dimension. For a fluid confined in a domain $\Omega_\varepsilon = I_\varepsilon \times (0,1)$, with $I_\varepsilon = (0,\varepsilon) \in \mathbb R$, the author analyzed the case where the viscosity coefficients have the following form
\begin{equation} \label{coeff-2d}
    0 < \mu = const., \ \ \lambda(\varrho) = b \varrho^\beta, 
\end{equation}
for $b>0$ and $\beta>3$, obtaining the convergence of the strong solution of the 2D system to the strong solution to the 1D system in the limit of $\varepsilon\rightarrow0$ through the use of a relative entropy inequality method. Existence results for strong solutions can be found in \cite{VK} and \cite{DWZ} respectively for 2D and 1D fluids. For other recent results concerning the dimension reduction analysis for compressible fluids the reader can refer, for example, to \cite{BFN}, \cite{BKM}, \cite{CDNS}, \cite{DCMP}, \cite{DNPB}, \cite{MN}, \cite{V}.

In our paper, we consider viscosity coefficients of the following form $\mu(\varrho)=\mu\varrho$, with $\mu > 0$ constant, and $\lambda(\varrho)=0$, and we extend the result of Zhang \cite{Z} showing that the weak solutions of the 3D system (\ref{cont}) - (\ref{mom}) converges to the strong solution of the respective 1D system as $\varepsilon\rightarrow0$.
Our result considers, in particular, the  so-called \textit{augmented version} of the Navier-Stokes system (\ref{cont}) - (\ref{mom}), see \cite{BDZ}, for which, recently, a new relative entropy inequality has been derived (see \cite{BNV}, \cite{BNV-2}). This new relative entropy inequality seems a suitable tools in order to perform a dimension reduction analysis for the density dependent viscosity case and, in general, for weak-strong uniqueness problems (see e.g. \cite{BGL}, \cite{CD}).

\section{Weak/$\kappa$-entropy solutions, augmented version of the Navier-Stokes system and main result}

In the following section we introduce the $\kappa$-entropy solutions and the augmented version of the compressible Navier-Stokes system, and we present the main result.

\subsection{\textit{$\kappa$-entropy} solutions}
In \cite{BDZ}, the authors introduced the so-called \textit{$\kappa$-entropy} solutions for the system (\ref{cont}) - (\ref{mom}).
\begin{defn} \label{ke-sol}
Let $T>0$ and $\kappa$ be such that $0<\kappa<1$. The couple $(\varrho,\vu)$ is called a global $\kappa$-entropy solution to the compressible Navier-Stokes system (\ref{cont}) - (\ref{mom}) if the following conditions are satisfied.

\bigskip
(1) Equations.

\bigskip
The continuity equation is satisfied in the following sense
\begin{equation} \label{mass-ke}
    -\int_0^T \int_{\Omega_\varepsilon} \varrho \partial_t \xi dxdt - -\int_0^T \int_{\Omega_\varepsilon} \varrho \vu \cdot \nabla \xi dxdt = \int_{\Omega_\varepsilon} \varrho(0) \xi(0) dx
\end{equation}
for all $\xi \in C_c^\infty([0,T)\times\Omega_\varepsilon)$.

\bigskip
The momentum equation is satisfied in the following sense
$$
-\int_0^T \int_{\Omega_\varepsilon} \varrho \vu \cdot \partial_t \phi dxdt -\int_0^T \int_{\Omega_\varepsilon} (\varrho \vu \otimes \vu) : \nabla \phi dxdt 
+\int_0^T \int_{\Omega_\varepsilon} 2\mu(\varrho) D(\vu) : \nabla \phi dxdt
$$ 
\begin{equation} \label{mom-ke}
    + \int_0^T \int_{\Omega_\varepsilon}
    \lambda (\varrho) \mbox{div} \vu \mbox{div} \phi dxdt
    - \int_0^T \int_{\Omega_\varepsilon} p(\varrho) \mbox{div} \phi dxdt
    = \int_{\Omega_\varepsilon} \varrho \vu(0) \cdot \phi(0) dx
\end{equation}
for all $\phi \in (C_c^\infty([0,T)\times\Omega_\varepsilon))^3$.

\bigskip
(2) Entropy estimate.

\bigskip
$(\varrho,\vu)$ satisfies, for all $t \in [0,T]$, the following $\kappa$-entropy estimate

$$
\sup_{t \in [0,T]} \left[ \int_{\Omega_\varepsilon}
\varrho \left( \frac{|\vu+2\kappa \nabla \varphi(\varrho)|^2}{2}+\frac{|2\sqrt{(1-\kappa)\kappa}\nabla\varphi|^2}{2}\right)(t)dx+\int_{\Omega_\varepsilon}\varrho e(\varrho)(t)dx \right]
$$

$$
+2 \kappa \int_0^T \int_{\Omega_\varepsilon} \mu(\varrho) |A(\vu)|^2 dxdt
+ 2 \kappa \int_0^T \int_{\Omega_\varepsilon} \frac{\mu'(\varrho)p'(\varrho)}{\varrho}|\nabla\varrho|^2 dxdt
$$

$$
2(1-\kappa) \int_0^T \int_{\Omega_\varepsilon} \mu(\varrho) |D(\vu)|^2 dxdt
+ 2(1-\kappa) \int_0^T \int_{\Omega_\varepsilon} (\mu'(\varrho)\varrho-\mu(\varrho))|\mbox{div}\vu|^2 dxdt
$$

\begin{equation} \label{ke-ineq}
    \leq
    \int_{\Omega_\varepsilon}
    \varrho \left( \frac{|\vu+2\kappa \nabla \varphi(\varrho)|^2}{2}+\frac{|2\sqrt{(1-\kappa)\kappa}\nabla\varphi|^2}{2}\right)(0)dx+\int_{\Omega_\varepsilon}\varrho e(\varrho)(0)dx
\end{equation}
with $\varphi'(\varrho)=\mu(\varrho)/\varrho$ and the internal energy $e(\varrho)$ defined by
$$
\frac{\varrho^2 de(\varrho)}{d\varrho}=p(\varrho).
$$
\end{defn}
\begin{rem} \label{BD}
The inequality (\ref{ke-ineq}) is a generalization of the well known BD-entropy introduced by Bresch and Desjardins in \cite{BD} (see also \cite{BD-1}) in the case $\kappa=1$. Moreover, it generalizes the one obtained in \cite{J} (see also \cite{GV}), taking $\kappa=1/2$ and $\varphi(\varrho)=\log \varrho$.   
\end{rem}

The existence of $\kappa$-entropy solutions for the compressible Navier-Stokes equations with $\mu(\varrho)=\mu\varrho$, with $\mu > 0$ constant, and $\lambda(\varrho)=0$ without extra terms (capillary, drag, singular pressure) have been recently proved in \cite{LX} and \cite{VY}. More precisely, the authors in \cite{LX} and \cite{VY} proved the global existence of weak solutions. However, as remarked in \cite{BDZ} (see Remark 2), a global weak solution of the system (\ref{cont}) - (\ref{mom}) is therefore a global $\kappa$-entropy solution. The converse is however not clear.

\subsection{Augmented Navier-Stokes system}

As observed in \cite{BDZ}, the system (\ref{cont}) - (\ref{mom}) can be reformulated through an augmented version. Indeed, for a coefficient $\kappa\in (0,1)$, defining the velocity $\vv = \vu + 2\kappa\mu \nabla\log\varrho$ and $\vw = 2\sqrt{\kappa(1-\kappa)}\mu \nabla \log \varrho$, the augmented version of the Navier-Stokes system reads as follows

\begin{equation} \label{cont-ag}
\pat \varrho + \diver (\varrho\vu) = 0,
\end{equation}

$$
\pat(\varrho \vv) + \diver (\varrho \vv \otimes \vu) + \nabla p(\varrho) 
$$
\begin{equation} \label{mom1-ag}
= \mu \diver (2\varrho (1-\kappa) D(\vv)) + \mu \diver (2\kappa \varrho A(\vv) - \mu\diver \left(2\sqrt{\kappa(1-\kappa)}\varrho \nabla \vw\right)
\end{equation}

$$
\pat(\varrho \vw) + \diver (\varrho \vw\otimes \vu ) 
$$
\begin{equation} \label{mom2-ag}
= \mu \diver (2\kappa\varrho \nabla \vw) - \mu \diver \left(2\sqrt{\kappa(1-\kappa)}\varrho (\nabla \vv^T)\right).
\end{equation}
As the authors remarked in \cite{BGL}, is important to mention that a global weak/$\kappa$-entropy solution $(\varrho, \vu)$ of the compressible Navier-Stokes system is also a solution of the augmented version.

\subsubsection{1D case}
In the one-dimensional case, we have
\begin{equation} \label{1d-cont}
\pat\overline\varrho + \pay(\overline\varrho \overline u)  = 0,
\end{equation}
\begin{equation} \label{1d-mom-1}
\pat(\overline\varrho \overline w) + \pay(\overline\varrho \overline w \overline u)  = -2\mu \pay(\sqrt{\kappa(1-\kappa)}\overline\varrho \pay \overline v) + \mu \pay(2\kappa \overline\varrho \pay \overline w),
\end{equation}
\begin{equation} \label{1d-mom-2}
\pat(\overline\varrho \overline v) + \pay (\overline \varrho \overline v \overline u) + \pay p(\overline  \varrho)  = (1-\kappa)2\mu \pay(\overline\varrho \overline v) - \sqrt{1-\kappa}\sqrt \kappa 2\mu \pay(\overline \varrho \pay \overline w),
\end{equation}
where all unknowns are overlined in order to distinguish them from their 3D counterparts. This system is considered on a time-space
$(0,T)\times (0,1)$ and, as in the previous case, it holds that
\begin{equation}
\overline v = \overline u + 2\kappa \mu \pay\log \overline\varrho,\qquad \overline w = 2\sqrt{\kappa(1-\kappa)}\mu \pay\log \overline\varrho.
\end{equation}
The assumed boundary conditions are of the form $\overline u(\cdot,0)= \overline u(\cdot,1)$ and $\overline\varrho(\cdot,0) = \overline\varrho(\cdot,1)$.

Hereinafter, we assume that the solution $(\overline\varrho,\overline u)$ to the 1D-case fulfills
\begin{equation}\label{eq:sol.1d}
\pay \overline u,\ \pay \log \overline \varrho,\ \pay^2 \log\overline\varrho \in L^\infty((0,T)\times(0,1)).
\end{equation}
The existence of solution which satisfies $\pay\overline u,\ \pay^2 \overline u,\ \pay \overline\varrho\in C^\alpha((0,T)\times (0,1))$ is known for the initial conditions satisfying $\overline u(0,\cdot) \in C^{1+\alpha}(0,1)$ and $\overline\varrho(0,\cdot)\in C^{2+\alpha}$. This was extensively studied -- we refer to the works by Kawohl \cite{K} or Kazhikhov and Shelukhin \cite{KaSh}. Although they worked with different boundary condition we do not expect any problem when adopting his method to the considered periodic boundary condition.
Nevertheless, we need to know further estimates concerning also $\pay^2 \overline\varrho$ -- this remains an open question.

\subsection{Relative entropy}

In the following we introduce a relative entropy inequality between a global weak/$\kappa-$entropy solution $(\varrho,\vv,\vw)$ to the system (\ref{cont-ag}) - (\ref{mom2-ag}) and another state of the fluid $(r,\vV,\vW)$ specified below. Let us note that the functions $(\varrho,\vv,\vw)$ also depend on $\varepsilon$, however this dependency is not emphasized in case we avoid any misunderstandings.

Inspired by \cite{BNV}, we define an entropy functional, $E$, for any triple of smooth functions $r,\vV,\vW:(0,T)\times \Omega_\varepsilon\mapsto (0,\infty)\times \mathbb R^3 \times \mathbb R^3$ by the following relation
$$
E(\varrho,\vv,\vw|r,\vV,\vW) 
$$
\begin{equation} \label{entr-funct}
    = \frac 12  \int_{\Omega_\varepsilon}\varrho(|\vw - \vW|^2 + |\vv - \vV|^2)\ {\rm d}x + \int_{\Omega_\varepsilon} (P(\varrho) - P(r) - P'(r)(\varrho - r))\ {\rm d}x
\end{equation}
where $P$ denotes a pressure potential given by 
\begin{equation} \label{P}
P(\varrho) = \varrho\int_1^\varrho \frac{p(s)}{s^2}\ {\rm d}s.
\end{equation}
Now, assuming that $r$, $r^{-1}$, $\vV$ and $\vW$ are continuously differentiable up to the boundary, is possible to derive the following entropy inequality (for the details of the derivation see \cite{BNV-2})
$$
E(\varrho,\vv,\vw|r,\vV,\vW)(\tau) - E(\varrho,\vv,\vw| r, \vV,\vW)(0)
$$
$$
+2k\mu \int_0^\tau \int_{\Omega_\varepsilon} \varrho |A(\vv - \vV)|^2\ {\rm d}x{\rm d}t 
+ 2\mu \int_0^\tau \int_{\Omega_\varepsilon} \varrho |D(\sqrt{(1-\kappa)} (\vv - \vV) - \sqrt\kappa(\vw - \vW)|^2\ {\rm d}x{\rm d}t
$$
$$
+2\kappa\mu \int_0^\tau \int_{\Omega_\varepsilon} \varrho \left[p'(\varrho)\nabla\log \varrho - p'(r)\nabla \log r\right] \cdot \left[\nabla \log \varrho - \nabla \log r\right]\ {\rm d}x{\rm d}t
$$
$$
\leq \int_0^\tau \int_{\Omega_\varepsilon} \varrho \left( \left(\vu \cdot \nabla \vW\right) \cdot \left(\vW - \vv\right)  + \left(\vu\cdot \nabla \vV\right)\cdot (\vV - \vv)\right)\ {\rm d}x{\rm d}t
$$
$$
+ \int_0^\tau \int_{\Omega_\varepsilon} \varrho \left(\pat \vW\cdot (\vW - \vw) + \pat \vV \cdot \left(\vV - \vv\right)\right)\ {\rm d}x{\rm d}t
$$
$$ 
+ \int_0^\tau \int_{\Omega_\varepsilon} \pat P'(r) (r - \varrho)\ {\rm d}x{\rm d}t - \int_0^\tau \int_{\Omega_\varepsilon} \nabla P'(r) \cdot \left[\varrho \vu - r \vU\right]\ {\rm d}x{\rm d}t
$$
$$ 
+ \int_0^\tau \int_{\Omega_\varepsilon} \left(p(r) - p(\varrho)\right)\diver \vU\ {\rm d}x{\rm d}t
$$
$$ 
- \kappa \int_0^\tau \int_{\Omega_\varepsilon} p'(\varrho)\nabla \varrho \cdot \left[2\mu \frac{\nabla r}{r} - \frac 1{\sqrt{(1-\kappa)\kappa}} \vW\right]\ {\rm d}x{\rm d}t
$$
$$ 
+ 2\mu \int_0^\tau \int_{\Omega_\varepsilon} \varrho \left( D(\sqrt{(1-\kappa)}\vV) - \nabla(\sqrt \kappa \vW) \right) : \left( D(\sqrt{(1-\kappa)} (\vV - \vv)) - \nabla (\sqrt \kappa (\vW - \vw))\right) \ {\rm d}x{\rm d}t
$$
$$
+ 2\kappa\mu \int_0^\tau \int_{\Omega_\varepsilon} \varrho A(\vV): A(\vV - \vv) \ {\rm d}x{\rm d}t + 2\kappa \mu \int_0^\tau \int_{\Omega_\varepsilon} \frac \varrho r p'(r)\nabla r\cdot \left(\frac{\nabla r}r - \frac{\nabla\varrho}{\varrho}\right)\ {\rm d}x{\rm d}t
$$
\begin{equation} \label{rei}
+2\sqrt{\kappa(1-\kappa)} \mu \int_0^\tau \int_{\Omega_\varepsilon} \varrho \left[ A(\vW): A(\vv - \vV) - A(\vw - \vW) : A(\vV)\right] \ {\rm d}x{\rm d}t
\end{equation}
for every $\tau \in (0,T)$, where we used the identities
\begin{equation}\label{uav}
\vu = \vv - \sqrt{\frac {\kappa}{1-\kappa}} \vw,\qquad \vU = \vV - \sqrt{\frac{\kappa}{1-\kappa}} \vW.
\end{equation}

\subsection{Main result}
Our main result reads.
\begin{thm} \label{thm}
Let $\overline\varrho,\ \overline v,\ \overline w$ be a classical solution to (\ref{1d-cont}) - (\ref{1d-mom-2})  emanating from the initial data $\overline\varrho_0,$ $\overline v_0$ and $\overline w_0$ which posses the regularity properties mentioned in \eqref{eq:sol.1d}. Further, let $\varrho_\varepsilon,$ $\vv_\varepsilon$, $\vw_\varepsilon$ be a global weak/$\kappa-$entropy solution to (\ref{cont-ag}) - (\ref{mom2-ag}) emanating from $\varrho_{0,\varepsilon},$ $\vv_{0,\varepsilon}$, $\vw_{0,\varepsilon}$. Let, moreover, the initial data satisfy
$$
\frac 1{|\Omega_\varepsilon|}  \int_{\Omega_\varepsilon}\frac 12 \varrho_{0,\varepsilon}(|\vw_{0,\varepsilon} - \overline w_0|^2 + |\vv_{0,\varepsilon} - \overline v_0|^2)\ {\rm d}x
$$
\begin{equation}
  + \frac 1{|\Omega_\varepsilon|} \int_{\Omega_\varepsilon} P(\varrho_{0,\varepsilon}) - P(\overline\varrho_0) - P'(\overline\varrho_0)(\varrho_{0,\varepsilon} - \overline\varrho_0)\ {\rm d}x \to 0
\end{equation}
Then 
\begin{equation}
(\varrho_\varepsilon, \vv_\varepsilon,\vw_\varepsilon)\to (\overline \varrho,\overline v,\overline w) \ \ as \ \ \varepsilon\rightarrow0
\end{equation}
in the following sense
\begin{equation}
\begin{split}
{\rm esssup}_{t\in(0,T)} \frac1{|\Omega_\varepsilon|}\|\varrho_\varepsilon - \overline\varrho\|_{\min\{2,\gamma\}}^{\min\{2,\gamma\}}& \to 0,\\
{\rm esssup}_{t\in(0,T)} \frac1{|\Omega_\varepsilon|}\int_{\Omega_\varepsilon} \varrho_\varepsilon\left(|\vw_\varepsilon - \overline w|^2 + |\vv_{\varepsilon} - \overline v|^2\right)\ {\rm d}x&\to 0.
\end{split}
\end{equation}
\end{thm}

\begin{rem}
We would like to mention that the convergence (2.9) does not imply $\vu_\varepsilon\rightarrow\overline{u}$, strictly speaking. However, we can redefine $\vw$, namely $\vw\rightarrow\vw\kappa/\sqrt{\kappa(1-\kappa)}$. Consequently, $\vu_\varepsilon\rightarrow\overline{u}=\overline{v}-\overline{w}$.
\end{rem}

\section{Proof of the Theorem \ref{thm}}

We set $x_3 = y$ and we use \eqref{rei} with
\begin{equation} \label{set}
\begin{split}
r(x_1,x_2,x_3) &= \overline\varrho(x_3),\\
\vV_1(x_1,x_2,x_3)  = 0, \ \ 
\vV_2(x_1,x_2,x_3)  = 0, \ \ 
\vV_3(x_1,x_2,x_3) & = \overline v(x_3),\\
\vW_1(x_1,x_2,x_3)  = 0, \ \ 
\vW_2(x_1,x_2,x_3)  = 0, \ \ 
\vW_3(x_1,x_2,x_3) & = \overline w(x_3).\\
\end{split}
\end{equation}
Consequently, \eqref{rei} may be rewritten as
$$
E(\varrho,\vv,\vw|r,\vV,\vW)(\tau) - E(\varrho,\vv,\vw| r, \vV,\vW)(0)\\
$$
$$
+2k\mu \int_0^\tau \int_{\Omega_\varepsilon} \varrho |A(\vv - \vV)|^2\ {\rm d}x{\rm d}t
+ 2\mu \int_0^\tau \int_{\Omega_\varepsilon} \varrho |D(\sqrt{(1-\kappa)} (\vv - \vV) - \sqrt\kappa(\vw - \vW)|^2\ {\rm d}x{\rm d}t\\
$$
$$
+2\kappa\mu \int_0^\tau \int_{\Omega_\varepsilon} \varrho \left[p'(\varrho)\nabla\log \varrho - p'(\overline\varrho)\nabla \log \overline\varrho\right] \cdot \left[\nabla \log \varrho - \nabla \log \overline\varrho\right]\ {\rm d}x{\rm d}t\\
$$
$$
\leq \int_0^\tau \int_{\Omega_\varepsilon} \varrho \left( \left(\vu_3 \pay \overline w\right) \left(\overline w - \vv_3\right)  + \left(\vu_3\pay \overline v\right) (\overline v - \vv_3)\right)\ {\rm d}x{\rm d}t\\
$$
$$
+ \int_0^\tau \int_{\Omega_\varepsilon} \varrho \left(\pat \overline w (\overline w - \vw_3) + \pat \overline v \cdot \left(\overline v - \vv_3\right)\right)\ {\rm d}x{\rm d}t\\
$$
$$
+ \int_0^\tau \int_{\Omega_\varepsilon} \pat P'(\overline\varrho) (\overline\varrho - \varrho)\ {\rm d}x{\rm d}t - \int_0^\tau \int_{\Omega_\varepsilon} \pay P'(\overline\varrho) \cdot \left[\varrho \vu_3 - \overline\varrho \overline u\right]\ {\rm d}x{\rm d}t\\
$$
$$
+ \int_0^\tau \int_{\Omega_\varepsilon} \left(p(\overline \varrho) - p(\varrho)\right)\pay \overline u\ {\rm d}x{\rm d}t\\
$$
$$
 - \kappa \int_0^\tau \int_{\Omega_\varepsilon} p'(\varrho)\nabla \varrho \cdot \left[2\mu \frac{\pay \overline \varrho}{\overline \varrho} - \frac 1{\sqrt{(1-\kappa)\kappa}} \overline w\right]\ {\rm d}x{\rm d}t\\
$$
$$
+ 2\mu \int_0^\tau \int_{\Omega_\varepsilon} \varrho \left( \pay(\sqrt{(1-\kappa)}\overline v - \sqrt \kappa \overline w) \right) \left( \pay(\sqrt{(1-\kappa)} (\overline v - \vv_3) - \sqrt \kappa (\overline w - \vw_3))\right) \ {\rm d}x{\rm d}t\\
$$
\begin{equation} \label{rei2}
+ 2\kappa \mu \int_0^\tau \int_{\Omega_\varepsilon} \frac \varrho {\overline\varrho} p'(\overline\varrho)\pay \overline \varrho \left(\frac{\pay \overline\varrho}{\overline\varrho} - \frac{\partial_{x_3}\varrho}{\varrho}\right)\ {\rm d}x {\rm d}t.
\end{equation}
We multiply \eqref{1d-mom-1} by $\frac{\varrho}{\overline \varrho}(\overline w - \vw_3)$ and \eqref{1d-mom-2} by $\frac\varrho{\overline\varrho}(\overline v - \vv_3)$ to deduce
\begin{multline} \label{1d-prod}
\int_{\Omega_\varepsilon}\left(\varrho \pat \overline w + \varrho \overline u \pay \overline w\right)(\overline w - \vw_3) + \left(\varrho \pat \overline v + \varrho \overline u \pay \overline v\right)(\overline v - \vv_3)\ {\rm d}x =  \\
= \int_{\Omega_\varepsilon} 2\mu \varrho\left(\sqrt{1-\kappa} \pay \overline v - \sqrt \kappa \pay \overline w\right) \left( \sqrt\kappa\pay (\overline w - \vw_3) - \sqrt{1-\kappa} (\overline v - \vv_3)\right)\ {\rm d}x\\
+ \int_{\Omega_\varepsilon} 2\mu \overline\varrho \pay \left(\frac \varrho {\overline\varrho}\right) \left[(\sqrt{1-\kappa} \pay\overline v - \sqrt \kappa \pay \overline w)(\sqrt \kappa (\overline w - \vw_3) - \sqrt{1-\kappa}(\overline v - \vv_3)\right]\ {\rm d}x\\
- \int_{\Omega_\varepsilon} \pay p(\overline\varrho) (\overline v - \vv_3)\frac{\varrho}{\overline\varrho}\ {\rm d}x.
\end{multline}
We plug this into \eqref{rei2} to get
\begin{multline} \label{rei3}
E(\varrho,\vv,\vw|r,\vV,\vW)(\tau) - E(\varrho,\vv,\vw| r, \vV,\vW)(0)\\
+2k\mu \int_0^\tau \int_{\Omega_\varepsilon} \varrho |A(\vv)|^2\ {\rm d}x{\rm d}t + 2\mu \int_0^\tau \int_{\Omega_\varepsilon} \varrho |D(\sqrt{(1-\kappa)} (\vv - \vV) - \sqrt\kappa(\vw - \vW)|^2\ {\rm d}x{\rm d}t\\
+2\kappa\mu \int_0^\tau \int_{\Omega_\varepsilon} \varrho \left[p'(\varrho)\nabla\log \varrho - p'(\overline\varrho)\nabla \log \overline\varrho\right] \cdot \left[\nabla \log \varrho - \nabla \log \overline\varrho\right]\ {\rm d}x{\rm d}t\\
\leq \int_0^\tau \int_{\Omega_\varepsilon} \varrho \left( \left((\vu_3 -\overline u ) \pay \overline w\right) \left(\overline w - \vv_3\right)  + \left((\vu_3 - \overline u)\pay \overline v\right) (\overline v - \vv_3)\right)\ {\rm d}x{\rm d}t\\
+ \int_0^\tau \int_{\Omega_\varepsilon} \pat P'(\overline\varrho) (\overline\varrho - \varrho)\ {\rm d}x{\rm d}t - \int_0^\tau \int_{\Omega_\varepsilon} \pay P'(\overline\varrho) \cdot \left[\varrho \vu_3 - \overline\varrho \overline u\right]\ {\rm d}x{\rm d}t\\
+ \int_0^\tau \int_{\Omega_\varepsilon} \left(p(\overline \varrho) - p(\varrho)\right)\pay \overline u\ {\rm d}x{\rm d}t\\
- \kappa \int_0^\tau \int_{\Omega_\varepsilon} p'(\varrho)\partial_{x_3} \varrho \cdot \left[2\mu \frac{\pay \overline \varrho}{\overline \varrho} - \frac 1{\sqrt{(1-\kappa)\kappa}} \overline w\right]\ {\rm d}x{\rm d}t\\
+ 2\kappa \mu \int_0^\tau \int_{\Omega_\varepsilon} \frac \varrho {\overline\varrho} p'(\overline\varrho)\pay \overline \varrho \left(\frac{\pay \overline\varrho}{\overline\varrho} - \frac{\partial_{x_3}\varrho}{\varrho}\right)\ {\rm d}x{\rm d}t\\
- \int_0^\tau \int_{\Omega_\varepsilon} \pay p(\overline\varrho) (\overline v - \vv_3) \frac\varrho{\overline\varrho}\ {\rm d}x{\rm d}t\\
+ \int_0^\tau \int_{\Omega_\varepsilon} 2\mu \overline\varrho \pay \left(\frac \varrho{\overline\varrho}\right) \left[(\sqrt{1-\kappa}\pay \overline v - \sqrt \kappa \pay \overline w )(\sqrt \kappa (\overline w - \vw_3) - \sqrt{1-\kappa}(\overline v - \vv_3)\right]\ {\rm d}x{\rm d}t.
\end{multline}
Since
\begin{equation}
|\vu_3 - \overline u|^2 = \left|\vv_3 - \sqrt{\frac{\kappa}{1-\kappa}} \vw_3 - \overline v + \sqrt{\frac{\kappa}{1-\kappa}} \overline w\right|^2\leq c \left(|\vv_3 - \overline v|^2 + |\vw_3 - \overline w|^2\right),
\end{equation}
we deduce with help of the H\"older inequality
\begin{multline}
\int_0^\tau \int_{\Omega_\varepsilon} \varrho \left( \left((\vu_3 -\overline u ) \pay \overline w\right) \left(\overline w - \vw_3\right)  + \left((\vu_3 - \overline u)\pay \overline v\right) (\overline v - \vv_3)\right)\ {\rm d}x\\
\leq c\left(\|\pay \overline w\|_\infty + \|\pay \overline v\|_\infty\right) \int_0^\tau E(\varrho,\vv,\vw|r,\vV,\vW)(t)\ {\rm d}t.
\end{multline}
By definition of $\overline w$ and $\vw$, we get
\begin{equation}
\pay\left(\frac\varrho{\overline\varrho}\right) = \frac{\varrho}{\overline\varrho}\left(\pay \log \varrho - \pay \log \overline\varrho\right) = \frac{\varrho}{\overline\varrho} \frac{1}{2\mu\sqrt\kappa\sqrt{1-\kappa}} \left(\vw_3 - \overline w\right).
\end{equation}
Consequently, the following estimate is obtained by a similar computation as the previous one.
\begin{multline}
\int_0^\tau \int_{\Omega_\varepsilon} 2\mu \overline\varrho \pay \left(\frac \varrho{\overline\varrho}\right) \\ \cdot \left[(\sqrt{1-\kappa}\pay \overline v - \sqrt \kappa \pay \overline w )(\sqrt \kappa (\overline w - \vw_3) - \sqrt{1-\kappa}(\overline v - \vv_3)\right]\ {\rm d}x{\rm d}t\\
\leq c(\|\pay \overline w\|_\infty + \|\pay \overline v\|_\infty)\int_0^\tau E(\varrho,\vv,\vw|r,\vV,\vW)\ {\rm d}t.
\end{multline}
Note that the second and the third term on the left hand side of \eqref{rei3} are positive and thus they can be neglected. As a result, \eqref{rei3} yields
$$
E(\varrho,\vv,\vw|r,\vV,\vW)(\tau) - E(\varrho,\vv,\vw| r, \vV,\vW)(0)\\
$$
$$
+2\kappa\mu \int_0^\tau \int_{\Omega_\varepsilon} \varrho \left[p'(\varrho)\nabla\log \varrho - p'(\overline\varrho)\nabla \log \overline\varrho\right] \cdot \left[\nabla \log \varrho - \nabla \log \overline\varrho\right]\ {\rm d}x{\rm d}t\\
$$
$$
\leq  \int_0^\tau \int_{\Omega_\varepsilon} \pat P'(\overline\varrho) (\overline\varrho - \varrho)\ {\rm d}x{\rm d}t - \int_0^\tau \int_{\Omega_\varepsilon} \pay P'(\overline\varrho) \cdot \left[\varrho \vu_3 - \overline\varrho \overline u\right]\ {\rm d}x{\rm d}t\\
$$
$$
+ \int_0^\tau \int_{\Omega_\varepsilon} \left(p(\overline \varrho) - p(\varrho)\right)\pay \overline u\ {\rm d}x{\rm d}t\\
$$
$$
- \kappa \int_0^\tau \int_{\Omega_\varepsilon} p'(\varrho)\partial_{x_3} \varrho \cdot \left[2\mu \frac{\pay \overline \varrho}{\overline \varrho} - \frac 1{\sqrt{(1-\kappa)\kappa}} \overline w\right]\ {\rm d}x{\rm d}t\\
$$
$$
+ 2\kappa \mu \int_0^\tau \int_{\Omega_\varepsilon} \frac \varrho {\overline\varrho} p'(\overline\varrho)\pay \overline \varrho \left(\frac{\pay \overline\varrho}{\overline\varrho} - \frac{\partial_{x_3}\varrho}{\varrho}\right)\ {\rm d}x{\rm d}t\\
$$
$$
- \int_0^\tau \int_{\Omega_\varepsilon} \pay p(\overline\varrho) (\overline v - \vv_3) \frac\varrho{\overline\varrho}\ {\rm d}x{\rm d}t\\
$$
\begin{equation} \label{rei4}
+ c\left(\|\pay \overline w\|_\infty + \|\pay  \overline v\|_\infty\right)\int_0^\tau E(\varrho,\vv,\vw|r,\vV,\vW)\ {\rm d}t.
\end{equation}
We will continue similarly as in \cite{BNV-2}. As $P''(s) = \frac{p'(s)}s$, we get
\begin{equation}
\pay p(\overline \varrho)(\overline v - \vv_3) \frac\varrho{\overline \varrho} = \varrho \pay \overline \varrho \frac{p'(\overline\varrho)}{\overline\varrho} (\overline v - \vv_3) = \varrho \pay P'(\overline\varrho)(\overline v - \vv_3).
\end{equation}
We multiply \eqref{1d-cont} by $P''(\overline\varrho)$ to get
\begin{equation}
\pat P'(\overline \varrho) + \pay P'(\overline\varrho) \overline u =  - P''(\overline\varrho)\overline\varrho \pay \overline u = -p'(\overline\varrho)\pay \overline u.
\end{equation}
We collect the previous two estimate and we use \eqref{uav} in order to deduce that
\begin{multline}
\int_0^\tau \int_{\Omega_\varepsilon}\pat P'(\overline\varrho)(\overline\varrho - \varrho) - \pay P'(\overline\varrho)(\varrho \vu_3 - \overline\varrho \overline u) - \pay p(\overline\varrho)(\overline v - \vv_3)\frac{\varrho}{\overline\varrho}\ {\rm d}x{\rm d}t\\
=\int_0^\tau \int_{\Omega_\varepsilon} \pat P'(\overline\varrho)(\overline\varrho - \varrho) + \pay P'(\overline\varrho)\overline u (\overline\varrho - \varrho)  + \pay P'(\overline\varrho)\varrho (\overline u - \vu_3) \\ - \varrho\pay P'(\overline\varrho)(\overline v - \vv_3)\ {\rm d}x{\rm d}t\\
= \int_0^\tau \int_{\Omega_\varepsilon} -p'(\overline\varrho) \pay \overline u (\overline\varrho - \varrho) - \frac{\sqrt\kappa}{\sqrt{1-\kappa}}\varrho \pay P'(\overline\varrho)(\overline w - \vw_3)\ {\rm d}x{\rm d}t.
\end{multline}
Consequently, \eqref{rei4} becomes
$$
E(\varrho,\vv,\vw|r,\vV,\vW)(\tau) - E(\varrho,\vv,\vw| r, \vV,\vW)(0)\\
$$
$$
+2\kappa\mu \int_0^\tau \int_{\Omega_\varepsilon} \varrho \left[p'(\varrho)\nabla\log \varrho - p'(\overline\varrho)\nabla \log \overline\varrho\right] \cdot \left[\nabla \log \varrho - \nabla \log \overline\varrho\right]\ {\rm d}x{\rm d}t\\
$$
$$
\leq  \int_0^\tau \int_{\Omega_\varepsilon} \left(p(\overline \varrho) - p(\varrho) - p'(\overline\varrho)(\overline\varrho - \varrho)\right)\pay \overline u\ {\rm d}x{\rm d}t\\
$$
$$
-\int_0^\tau\int_{\Omega_\varepsilon} \frac{\sqrt \kappa}{\sqrt{1-\kappa}} \varrho\pay P'(\overline\varrho)(\overline w - \vw_3)\ {\rm d}x{\rm d}t\\
$$
$$
- \kappa \int_0^\tau \int_{\Omega_\varepsilon} p'(\varrho)\partial_{x_3} \varrho \cdot \left[2\mu \frac{\pay \overline \varrho}{\overline \varrho} - \frac 1{\sqrt{(1-\kappa)\kappa}} \overline w\right]\ {\rm d}x{\rm d}t\\
$$
$$
+ 2\kappa \mu \int_0^\tau \int_{\Omega_\varepsilon} \frac \varrho {\overline\varrho} p'(\overline\varrho)\pay \overline \varrho \left(\frac{\pay \overline\varrho}{\overline\varrho} - \frac{\partial_{x_3}\varrho}{\varrho}\right)\ {\rm d}x{\rm d}t\\
$$
\begin{equation} \label{rei5}
+ c\left(\|\pay \overline w\|_\infty + \|\pay  \overline v\|_\infty\right)\int_0^\tau E(\varrho,\vv,\vw|r,\vV,\vW)\ {\rm d}t.
\end{equation}
Since (see \cite[Section 4]{FJN})
\begin{equation}\label{press.est}
P(\varrho) - P(\overline\varrho) - P'(\overline\varrho)(\varrho - \overline\varrho) \sim\left\{
\begin{array}{l}
(\varrho - \overline\varrho)\ \mbox{for } \frac 12 \min (\overline\varrho)\leq \varrho \leq 2\max (\overline\varrho),\\
(1+\varrho^\gamma)\ \mbox{otherwise}.
\end{array}
\right.
\end{equation}
we get
\begin{equation}
p(\varrho) - p(\overline\varrho) - p'(\overline\varrho)(\varrho -\overline \varrho)\sim P(\varrho) - P(\overline\varrho) - P'(\overline\varrho)(\varrho - \overline\varrho)
\end{equation}
and consequently
\begin{equation}
\int_{\Omega_\varepsilon} (p(\overline\varrho) - p(\varrho) - p'(\overline\varrho)(\overline\varrho - \varrho))\pay \overline u\ {\rm d}x \leq c\|\pay \overline u\|_\infty E(\varrho,\vv,\vw|r,\vV,\vW).
\end{equation}
Due to the definition of $\overline w$, we have
\begin{equation}
2\mu \frac{\pay \overline\varrho}{\overline\varrho} - \frac 1{\sqrt{(1-\kappa)\kappa}} \overline w = 0.
\end{equation}
Further, we have
\begin{multline}
\frac{\sqrt\kappa}{\sqrt{1-\kappa}} \varrho \pay P'(\overline\varrho) (\overline w - \vw_3) \\ = 
\frac{\sqrt \kappa}{\sqrt{1-\kappa}} \varrho p'(\overline\varrho)\pay \overline\varrho \frac 1{\overline\varrho}\left(2\sqrt{\kappa(1-\kappa)} \mu \left(\pay \log \overline\varrho - \pay \log \varrho\right)\right)\\
\\
= 2\kappa\mu \frac{\varrho}{\overline\varrho}p'(\overline\varrho)\pay \overline\varrho (\pay \log\overline\varrho - \pay \log \varrho).
\end{multline}
Consequently, the sum of the second and fourth term on the right hand side of (\ref{rei5}) is zero.
It remains to deal with the second term on the left hand side. Here, we use the fact that $p(\varrho) = \varrho^\gamma$ (assuming $a=1$ without loss of generality). We have, 
\begin{multline}
\varrho (p'(\varrho)\pay \log \varrho - p'(\overline \varrho)\pay \log \varrho)(\pay \log \varrho - \pay \log \overline\varrho) \\
= \varrho p'(\varrho) |\pay \log \varrho - \pay \log \overline\varrho|^2 + \varrho (p'(\varrho) - p'(\overline\varrho))\pay \log \overline\varrho (\pay \log \varrho - \pay \log \overline\varrho)\\
=\varrho p'(\varrho) |\pay \log \varrho - \pay \log \overline\varrho|^2 + \pay(p(\varrho) - p(\overline\varrho) - p'(\overline\varrho)(\varrho - \overline\varrho)) \pay \log \overline\varrho\\ - (\varrho(p'(\varrho) - p'(\overline\varrho)) - p''(\overline\varrho)(\varrho - \overline\varrho)\overline\varrho)|\pay\log\overline\varrho|^2.
\end{multline}
The first term on the right hand side is definitely positive and thus it can be neglected in \eqref{rei5}. For the rest, by integration by parts, we have
\begin{multline}
\int_{\Omega_\varepsilon}\pay(p(\varrho) - p(\overline\varrho) - p'(\overline\varrho)(\varrho - \overline\varrho)) \pay \log \overline\varrho - (\varrho(p'(\varrho) - p'(\overline\varrho)) \\
- p''(\overline\varrho)(\varrho - \overline\varrho)\overline\varrho)|\pay\log\overline\varrho|^2\ {\rm d}x \\
 \leq c \int_{\Omega_\varepsilon} |(p(\varrho) - p(\overline \varrho) - p'(\overline\varrho) (\varrho - \overline\varrho)| |\pay^2\log \overline\varrho|\ {\rm d}x\\
 + c\int_{\Omega_\varepsilon} |\varrho (p'(\varrho) - p'(\overline\varrho)) - p''(\overline\varrho)(\varrho - \overline\varrho)\overline\varrho| |\pay \log \overline\varrho|^2 \ {\rm d}x.
\end{multline}
We end up with
\begin{multline}
E(\varrho,\vv,\vw|r,\vV,\vW)(\tau) - E(\varrho,\vv,\vw|r,\vV,\vW)(0) \\ \leq c(\|\pay \overline w\|_\infty,\|\pay \overline v\|_\infty, \|\pay \log \overline\varrho\|_\infty, \|\pay^2 \log \overline\varrho\|_\infty) \int_0^\tau E(\varrho,\vv,\vw|r,\vV,\vW)(t)\ {\rm d}t.
\end{multline}
Now the dependency on $\varepsilon$ matters and thus we use a index $\varepsilon$ to denote it. By the Gronwall inequality and dividing by $1/|\Omega_\varepsilon|$, we deduce that
$$
\frac{1}{|\Omega_\varepsilon|}E(\varrho_\varepsilon,\vv_\varepsilon,\vw_\varepsilon|r,\vV,\vW)(\tau) 
$$
\begin{equation}\label{eq:after.gronwall}
\leq \frac{c}{|\Omega_\varepsilon|}(\|\pay \overline w\|_\infty, \|\pay \overline v\|_\infty, \|\pay \log \overline \varrho\|_\infty, \|\pay^2 \log \overline\varrho\|_\infty, T) E(\varrho_\varepsilon,\vv_\varepsilon,\vw_\varepsilon|r,\vV,\vW)(0).
\end{equation}
for every $\tau\in [0,T]$ and for some constant $c$ independent of $\varepsilon$. Due to assumptions, we have
\begin{equation}
\frac{1}{|\Omega_\varepsilon|} E(\varrho_\varepsilon,\vv_\varepsilon,\vw_\varepsilon|r,\vV,\vW)(0) \rightarrow 0, \ \ \mbox{as} \ \ \varepsilon\rightarrow0. 
\end{equation}
Consequently, \eqref{eq:after.gronwall} yields the demanded convergencies and Theorem \ref{thm} is proved.

\bigskip
\section*{Acknowledgements}
M. C. has been supported by the Croatian Science Foundation under the project MultiFM IP-2019-04-1140. V. M. claims support of Czech Science Foundation project number GA19-04243S and RVO: 67985840

\end{document}